\documentclass[a4paper]{article}
\usepackage{mathbbold}
\usepackage{amsfonts}
\usepackage{amssymb}
\usepackage{mathrsfs}
\usepackage{float}
\usepackage{amsmath}
\usepackage{amsfonts,amssymb}
\usepackage{graphicx,color}
\usepackage{psfig}
\usepackage{mathrsfs}

\newtheorem{Theorem}{Theorem}[section]
\newtheorem{Definition}{Definition}[section]
\newtheorem{Lemma}{Lemma}[section]
\newtheorem{Example}{Example}[section]

\newtheorem{Remark}{Remark}[section]

\newtheorem{Proposition}{Proposition}[section]

\catcode`\@=11 \@addtoreset{equation}{section}

\catcode`\@=12

\def\2{{I \hskip -1.0mm I}}
\def\3{{I \hskip -1.0mm I\hskip -1.0mm I}}
\def\4{{I \hskip -0.9mm V}}
\def\6{{V \hskip -1.35mm I}}

\date{ }
\title{Hyperbolic mean curvature flow: Evolution of plane curves}
\author{De-Xing Kong\footnote{Center of Mathematical Sciences, Zhejiang University,
Hangzhou 310027, China;},$\quad$ Kefeng Liu\footnote{ Department of
Mathematics, University of California at Los Angeles, CA 90095,
USA;} $\quad$ and $\quad$  Zeng-Gui Wang\footnote{Department of
Mathematics, Shanghai Jiao Tong University, Shanghai 200240, China.}
\footnote{Corresponding author.}}
\date{ }
\begin{document}
\maketitle
\begin{abstract}
In this paper we investigate the one-dimensional hyperbolic mean
curvature flow for closed plane curves. More precisely, we consider
a family of closed curves $F:S^{1}\times [0,T)\to\mathbb{R}^{2}$
which satisfies the following evolution equation
$$\dfrac{\partial^{2} F}{\partial t^{2}}(u,t)
=k(u,t)\vec{N}(u,t)-\triangledown \rho(u,t),\ \quad \forall\
(u,t)\in S^{1}\times[0,T)$$  with the initial data
$$F(u,0)=F_{0}(u) \quad {\rm{ and}}\quad
\dfrac{\partial F}{\partial t}(u,0)=f(u)\vec{N_{0}},$$ where $k$ is
the mean curvature and $\vec{N}$ is the unit inner normal vector of
the plane curve $F(u,t)$, $f(u)$ and $\vec{N_{0}}$ are the initial
velocity and the unit inner normal vector of the initial convex
closed curve $F_{0}$ respectively, and $\triangledown\rho$ is given
by
$$\triangledown \rho\triangleq\left\langle\dfrac{\partial^{2} F}{\partial s\partial t},
\dfrac{\partial F}{\partial t}\right\rangle\vec{T},$$ in which
$\vec{T}$ stands for the unit tangent vector. The above problem is
an initial value problem for a system of partial differential
equations for $F$, it can be completely reduced to an initial value
problem for a single partial differential equation for its support
function. The latter equation is a hyperbolic Monge-Amp\`{e}re
equation. Based on this, we show that there exists a class of
initial velocities such that the solution of the above initial value
problem exists only at a finite time interval $[0,T_{\max })$ and
when $t$ goes to $T_{\max}$, the solution converges to a point.
%Furthermore, we also consider the hyperbolic mean curvature flow
%with the dissipative terms and prove the global existence of smooth
%solutions with small initial data.
In the end, we discuss the close
relationship between the hyperbolic mean curvature flow and the
equations for the evolving relativistic string in the Minkowski
space-time $\mathbb{R}^{1,1}$.

\vskip 6mm

\noindent{\bf Key words and phrases}: hyperbolic mean curvature
flow, hyperbolic Monge-Amp\`{e}re equation, closed plane curve,
short-time existence.

\vskip 3mm

\noindent{\bf 2000 Mathematics Subject Classification}: 58J45,
58J47.
\end{abstract}
\newpage
\baselineskip=7mm

\section{Introduction}

In this paper we study the closed convex evolving plane curves. More
precisely, we consider the following initial value problem
\begin{equation}\label{1.6}\left\{\aligned&\dfrac{\partial^{2} F}{\partial t^{2}}(u,t)
=k(u,t)\vec{N}(u,t)-\triangledown \rho(u,t),\quad \forall \ (u,t)\in S^{1}\times[0,T),\\
&F(u,0)=F_{0}(u),\\
&\dfrac{\partial F}{\partial
t}(u,0)=f(u)\vec{N_0},\endaligned\right.\end{equation} where $k$ is
the mean curvature, $\vec{N}$ is the unit inner normal at $F(u,t)$,
$F_{0}$ stands for a smooth strictly convex closed curve, $f(u)(\geq
0)$ and $\vec{N_{0}}$ are the initial velocity and inner normal
vector of $F_{0}$, respectively, and with $\vec{T}$ denoting the
unit tangent vector and $s$ the arclength parameter, $\triangledown
\rho$ is defined by
\begin{equation}\label{1.7}\triangledown \rho\triangleq\left\langle\dfrac{\partial^{2}
F}{\partial s\partial t}, \dfrac{\partial F}{\partial
t}\right\rangle\vec{T}.\end{equation}

This system is an initial value problem for a system of partial
differential equations for $F$, which can be completely reduced to
an initial value problem for a single partial differential equation
for its support function. The latter equation is a hyperbolic
Monge-Amp\`{e}re equation. Our first result is the following local
existence theorem for the initial value problem (\ref{1.6}).

\begin{Theorem}\label{2.1} {\bf (Local existences and uniqueness)}
Suppose that $F_{0}$ is a smooth strictly convex closed curve. Then
there exist a positive $T$ and a family of strictly convex closed
curves $F(\cdot,t)$ with $t\in [0,T)$ such that $F(\cdot,t)$
satisfies (\ref{1.6}), provided that $f(u)$ is a smooth function on
$S^{1}$.
\end{Theorem}

Our second main result is the following theorem.

\begin{Theorem}\label{4.1} Suppose that $F_{0}$ is a smooth strictly
convex closed curve. Then there exists a class of the initial
velocities such that the solution of (\ref{1.6}) with $F_{0}$ and
$f$ as initial curve and initial velocity of the initial curve,
respectively, exists only at a finite time interval $[0,T_{\max})$.
Moreover, when $t\to T_{\max}$, the solution $F(\cdot,t)$ converges
to a point.\end{Theorem}

After rescaling as Gage and Hamilton did, we can see that the
limiting solution will be a circle. We will also introduce
hyperbolic mean curvature flow with dissipative terms. A close
relation between our hyperbolic mean curvature flow and the string
evolving in the Minkowski space-time will be derived in the last
section of the paper.

For reader's convenience, we briefly discuss some history of
parabolic and hyperbolic mean curvature flows.

The parabolic theory for the evolving of plane curves, which in its
simplest form is based on the curve shortening equation
\begin{equation}\label{1.1} v=k\end{equation}
relating the normal velocity $v$ and the curvature $k$, has been
extremely successful in providing geometers with great insight. For
example, Gage and Halmilton \cite{gh} proved that, when the
 curve is strictly convex, the deformation decreases the isoperimetric
 ratio, and furthermore if it shrinks to a point $p$, the
 normalized curves, obtained by ``blowing up" the curves at $p$ so
 that its enclosed areas is $\pi$, must tend to the unit circle in a
 certain sense. Grayson \cite{g} generalized this result and showed that a
 smooth embedded plane curve first becomes convex and then shrinks
 to a point in a finite time. These results can be applied to many physical problems
 such as crystal growth, computer vision and image processing. One of the important applications
of mean curvature flow is that Huisken and Ilmanen developed a
theory of weak solutions of the inverse mean curvature flow and used
it to prove successfully the Riemannian Penrose inequality  which
plays an important role in general relativity (see \cite{hu2}).

However, to our knowledge, there is very few hyperbolic versions of
mean curvature flow. Melting crystals of helium exhibits a
phenomenon generally not found in other materials: oscillations of
the solid-liquid interface in which atoms of the solid move only
when they melt and enter the liquid (see \cite{g} and references
therein). Gurtin and Podio-Guidugli \cite{g} developed a hyperbolic
theory for the evolution of plane curves. Rostein, Brandon and
Novick-Cohen \cite{R} studied a hyperbolic theory by the mean
curvature flow equation \begin{equation}\label{1.2}v_{t}+\psi
v=k,\end{equation}where $v_{t}$ is the normal acceleration of the
interface, $\psi$ is a constant. A crystalline algorithm was
developed for the motion of closed polygonal curves.

The hyperbolic version of mean curvature flow is important in both
mathematics and applications, and has attracted many mathematicians
 to study it. He, Kong and Liu \cite{h} introduced hyperbolic mean
curvature flow from geometric point of view. Let $\mathscr{M}$ be a
Riemannian manifold and $X(\cdot,t):\mathscr{M}\rightarrow
\mathbb{R}^{n+1}$ be a smooth map. When $X$ is an isometric
immersion, the Laplacian of $F$ is given by $\triangle X=H\vec{N}$,
where $H$ is the mean curvature (i.e., the trace of the second
fundamental form) and $\vec{N}$ is the unit inner  normal vector.
The hyperbolic mean curvature flow is the following partial
differential equation of second order
\begin{equation}\label{1.3}\dfrac{\partial^2 }{\partial
t^2}X(u,t)=\triangle X \quad {\rm or} \quad \dfrac{\partial^2
}{\partial t^2}X(u,t)=H(u,t) \vec{N}(u,t),\quad \forall ~u\in
\mathscr{M} ,\quad
 \forall ~t>0.\end{equation}
 $X=X(u,t)$ is called a solution of
the hyperbolic mean curvature flow if it satisfies the equation
(\ref{1.3}). He, Kong and Liu in \cite{h} proved that the
corresponding system of partial differential equations are strictly
hyperbolic, and based on this, they also showed that this flow
admits a unique short-time smooth solution and possesses the
nonlinear stability defined on the Euclidean space with dimension
larger than 4. Moreover, the nonlinear wave equations satisfied by
curvatures are also derived in \cite{h}, these equations will play
an important role in future study. The hyperbolic mean curvature
flow was considered as one of the general hyperbolic geometric flows
introduced by Kong and Liu, see \cite{KL} for more discussions for
related hyperbolic flows and their applications to geometry and
Einstein equations.

Recently, Lefloch and Smoczyk \cite{l} studied the following
geometric evolution equation of hyperbolic type which governs the
evolution of a hypersurface moving in the direction of its mean
curvature vector
\begin{equation}\label{1.4}\left\{\aligned&\dfrac{\partial^2 }{\partial
t^2}X=eH(u,t)\vec{N}-\bigtriangledown e,\\& X(u,0)=X_{0},
\\
&\left(\dfrac{\partial X}{\partial
t}\right)_{t=0}^{\vec{T}_{0}}=0,\endaligned \right.\qquad
\end{equation}
where $\vec{T}_{0}$ stands for the unit tangential vector of the
initial hypersurface $X_0$,
$e\triangleq\dfrac{1}{2}\left(\left|\dfrac{d}{dt}X\right|^{2}+n\right)$
is the local energy density and $\bigtriangledown
e\triangleq\bigtriangledown^{i}e_{i}$, in which
$e_{i}=\dfrac{\partial e}{\partial x^{i}}$.
 This flow stems from a geometrically natural action containing
kinetic and internal energy terms. They have shown that the normal
hyperbolic mean curvature flow will blow up in finite time. In the
case of graphs, they introduce a concept of weak solution suitably
restricted by an entropy inequality and proved that the classical
solution is unique in the larger class of entropy solutions. In the
special case of one-dimensional graphs, a global-in-time existence
result is established. Moreover, an existence theorem has been
established under the assumption that the BV norm of initial data is
small.

The paper is organized as follows. In Section 2, a hyperbolic
Monge-Amp\`{e}re equation will be derived and a theorem on local
existence and uniqueness of the solution, i.e., Theorem 1.1 will be
proved. In Section 3, an example is given and then some properties
of the evolving curve have been established. The main result
--- Theorem 1.2 will be proved in Section 4. In section 5, we
consider the normal hyperbolic mean curvature flow with the
dissipative term and get the hyperbolic equation of $S$ and $k$
respectively. Section 6 is devoted to illustrating the relations
between the hyperbolic mean curvature flow and the evolution
equations for the relativistic string in the Minkowski space
$\mathbb{R}^{1,1}$.

\section{Hyperbolic Monge-Amp\`{e}re equation}

Roughly speaking, an evolving curve is a smooth family of curves
$u\mapsto F(u,t)$, where $u\in S^{1}$ and $t\in [0, T)$, in which
$T$ is called the duration of $F$. For a given curve $F(\cdot, t)$,
the underlying physics must be independent of the choice of the
parameter $u$, and hence can involve $F$ only through intrinsic
quantities such as curvature, normal acceleration and normal
velocity, which are independent of parametrization. On the other
hand, this invariance allows us to use any convenient
parametrization. The following notion is needed in our study.
\begin{Definition}\label{1.1} A curve $F: S^{1}\times [0,T)\to
\mathbb{R}^{2}$ evolves normally if
\begin{equation}\label{1.5}\left\langle\dfrac{\partial F}{\partial t},\dfrac{\partial
F}{\partial u}\right\rangle=0\end{equation} for all $(u,t)\in
S^{1}\times [0,T)$.
\end{Definition}

Definition 2.1 can be found in \cite {A} and \cite{l}. In this
paper, we shall restrict our attention to the parametrization
(\ref{1.5}). Such a parametrization significantly simplifies the
analysis. The normally evolving curve was first investigated by
Angenent and Gurtin \cite{A} and then further studied by Lefloch and
Smoczyk \cite{l}. The following result is important, it shows that
within a large class of time-dependent curves there is no essential
loss of generality in limiting attention to curves that evolves
normally.

\vskip2mm

\noindent{\bf Lemma A} {\it If the evolving curve $\mathscr{C}$ is
closed, then there is a parameter change $\phi$ for $\mathscr{C}$
such that $\mathscr{C}\circ\phi$ is a normally evolving curve.}

\vskip2mm

This important lemma was proved in \cite{A}. In fact, for the
initial value problem (\ref{1.6}), the initial velocity field is
normal to the curve, it can be proved that this property is
preserved during the evolution, that is to say, the flow (1.1) is
automatically a normal flow\footnote{Although the partial
differential equation in (2.2) only contains the first order
derivative of $F$ with respect to $t$, i.e., $F_t$, it is non-local
partial differential equation, and it is not easier to handle than
the second order partial differential equation in (1.1).}
\begin{equation}\label{1.8}\left\{\aligned&\dfrac{\partial F}{\partial
t}=\sigma(u,t)\vec{N},\\
&F(u,0)=F_{0}(u),\endaligned\right.\end{equation} where
$\sigma(u,t)=f(u)+\int_0^{t}k(u,\xi)d\xi$, hence we have
\begin{equation}\label{1.9}\dfrac{\partial\sigma}{\partial
t}=k(u,t),\quad \sigma\dfrac{\partial\sigma}{\partial
s}=\left\langle\dfrac{\partial^{2} F}{\partial s\partial t},
\dfrac{\partial F}{\partial t}\right\rangle.
\end{equation}
Here we denote by $s=s(\cdot,t)$ the arclength parameter of the
curve $F(\cdot,t): S^{1}\to \mathbb{R}^{2}$. The operator
$\partial/\partial s$ is given in terms of $u$ by
$$\dfrac{\partial}{\partial s}=\frac{1}{\upsilon}\dfrac{\partial}{\partial u},$$
where $$\upsilon=\sqrt{(\partial x/\partial u)^{2}+(\partial
y/\partial u)^{2}}=|\partial F/\partial u|.$$ By Frenet formula,
$$\dfrac{\partial\vec{T}}{\partial s}=k\vec{N},\quad \dfrac{\partial\vec{N}}
{\partial s}=-k\vec{T}.$$
 Then $\{\vec{T},\ \vec{N}\}$ is an orthogonal basis of $\mathbb{R}^{2}$. Let us denote $\theta$
 to be the unit outer normal angle for a convex closed curve $F\ :S^{1}\to
 \mathbb{R}^{2}$. Hence, $$\vec{N}=(-\cos \theta,-\sin \theta),\quad \vec{T}=(-\sin \theta, \cos \theta),$$
 and by Frenet
 formula, we have $$\dfrac{\partial\theta}{\partial s}=k.$$
 Furthermore,
 $$\dfrac{\partial N}{\partial t}=-\dfrac{\partial \theta}{\partial t}\vec{T},\quad \
\dfrac{\partial T}{\partial t}=\dfrac{\partial \theta}{\partial
t}\vec{N}. $$ Using the previous definition, we have
$$\dfrac{\partial^{2}}{\partial t\partial s}=k \sigma\dfrac{\partial}{\partial
s}+\dfrac{\partial^{2}}{\partial s\partial t},$$ and observing
$$\vec{T}(u,t)=\dfrac{\partial F}{\partial s}(u,t),$$ we deduce
$$\dfrac{\partial\vec{T}}{\partial t}=\dfrac{\partial }{\partial t}\left(\dfrac{\partial F}{\partial s}\right)
=\dfrac{\partial\sigma}{\partial s}\vec{N}$$ and
$$\dfrac{\partial\vec{N}}{\partial t}=-\dfrac{\partial\sigma}{\partial s}\vec{T},$$
hence,
$$\dfrac{\partial\theta}{\partial t}=\dfrac{\partial\sigma}{\partial s}.$$

Suppose that $F(u,t):\;S^{1}\times[0,T)\rightarrow \mathbb{R}^{2}$
is a family of convex curves satisfying the curve shortening flow
$(\ref{1.6})$. Let us use the normal angle to parameterize each
convex curve $F(\cdot,t)$, i.e., set
$$\widetilde{F}(\theta,\tau)=F(u(\theta,\tau), t(\theta,\tau)),$$
where $t(\theta,\tau)=\tau$. Here, $\vec{N}$ and $\vec{T}$ are
independent of the parameter $\tau$, which can be proved as follows:
$$0=\dfrac{\partial \theta}{\partial\tau}=\dfrac{\partial\theta}{\partial u}
\dfrac{\partial u}{\partial \tau}+\dfrac{\partial\theta}{\partial
t},$$ then $$\dfrac{\partial\theta}{\partial t}=-\dfrac{\partial
\theta}{\partial u}\dfrac{\partial u}{\partial
\tau}=-\dfrac{\partial \theta}{\partial s}\dfrac{\partial
s}{\partial u}\dfrac{\partial u}{\partial
\tau}=-k\upsilon\dfrac{\partial u}{\partial \tau}.$$ Hence,
$$\dfrac{\partial\vec{T}}{\partial\tau}=\dfrac{\partial\vec{T}}{\partial t}+
\dfrac{\partial\vec{T}}{\partial u}\dfrac{\partial u}{\partial
\tau}=\dfrac{\partial\theta}{\partial
t}\vec{N}+\dfrac{\partial\vec{T}}{\partial s}\dfrac{\partial
s}{\partial u}\dfrac{\partial u}{\partial
\tau}=\left(\dfrac{\partial\theta}{\partial
t}+k\upsilon\dfrac{\partial u}{\partial \tau}\right)\vec{N}=0 $$ and
$$\dfrac{\partial\vec{N}}{\partial\tau}=\dfrac{\partial\vec{N}}{\partial t}+
\dfrac{\partial\vec{N}}{\partial u}\dfrac{\partial u}{\partial
\tau}=-\dfrac{\partial\theta}{\partial
t}\vec{T}+\dfrac{\partial\vec{N}}{\partial s}\dfrac{\partial
s}{\partial u}\dfrac{\partial u}{\partial
\tau}=-\left(\dfrac{\partial\theta}{\partial
t}+k\upsilon\dfrac{\partial u}{\partial \tau}\right)\vec{T}=0.$$
 By
the chain rule,
$$\dfrac{\partial\widetilde{F}}{\partial \tau}=\dfrac{\partial F}
{\partial u}\dfrac{\partial u}{\partial\tau}+\dfrac{\partial
F}{\partial t}$$ and
$$\dfrac{\partial^{2}\widetilde{F}}{\partial \tau^{2}}=\dfrac{\partial F}
{\partial u}\dfrac{\partial^{2} u}{\partial
\tau^{2}}+\dfrac{\partial^{2} F} {\partial
u^{2}}\left(\dfrac{\partial
u}{\partial\tau}\right)^{2}+2\dfrac{\partial^{2}F}{\partial
u\partial t}\dfrac{\partial u}{\partial \tau}+\dfrac{\partial^{2}
F}{\partial t^{2}}.$$

 The support function of $F$ is given by
 $$\aligned S(\theta,\tau)=&\langle\widetilde{F}(\theta,\tau),-\vec{N}\rangle,\\
 =&\langle \widetilde{F}(\theta,\tau),(\cos\theta,\sin\theta)\rangle\\
 =&x(\theta,\tau)\cos\theta+y(\theta,\tau)\sin\theta.\endaligned$$
Its derivative satisfies
$$\aligned
S_{\theta}(\theta,\tau)=&-x(\theta,\tau)\sin\theta+y(\theta,\tau)\cos\theta+
\langle
\widetilde{F}_{\theta}(\theta,\tau),(\cos\theta,\sin\theta)\rangle\\
=&-x(\theta,\tau)\sin\theta+y(\theta,\tau)\cos\theta\\
=&\langle\tilde{F}(\theta,\tau),\vec{T}\rangle,\endaligned$$ where
$$\langle
\widetilde{F}_{\theta}(\theta,\tau),(\cos\theta,\sin\theta)\rangle=0,$$
namely, the tangent vector is orthogonal to the unit normal vector.
And then the curve can be represented by the support function
\begin{equation}\left\{\aligned x=&S\cos\theta-S_{\theta}\sin\theta,
\\y=&S\sin\theta+S_{\theta}\cos\theta.\endaligned\right.\end{equation}
Thus all geometric quantities of the curve can be represented by the
support function. In particular, the curvature can be written as
$$k=\dfrac{1}{S_{\theta\theta}+S}.$$
In fact, by the definition of the support function,
$$\aligned S_{\theta\theta}+S=&-x_{\theta}\sin\theta+y_{\theta}\cos\theta
-x\cos\theta-y\sin\theta+x\cos\theta+y\sin\theta\\
=&\left\langle\dfrac{\partial
\widetilde{F}}{\partial\theta},\vec{T}\right\rangle=\left\langle\dfrac{\partial
\widetilde{F}}{\partial s}\dfrac{\partial s}{\partial
\theta},\vec{T}\right\rangle=\dfrac{1}{k}.\endaligned$$

 We know that the support function
$$S(\theta,\tau)=\left\langle\widetilde{F}(\theta,\tau),-\vec{N}\right\rangle$$
satisfies
$$\aligned S_{\tau}=&\left\langle\dfrac{\partial\widetilde{F}}{\partial
\tau},-\vec{N}\right\rangle+\left\langle\widetilde{F}(\theta,\tau),\dfrac{\partial\vec{N}}{\partial\tau}\right\rangle\\
=&\left\langle\dfrac{\partial F}{\partial u}\dfrac{\partial
u}{\partial \tau}+\dfrac{\partial F}{\partial
t},-\vec{N}\right\rangle\\=&\langle\dfrac{\partial F}{\partial
t},-\vec{N}\rangle=-\tilde{\sigma}(\theta,\tau),\endaligned$$ where
$$\left\langle\widetilde{F}(\theta,\tau),\dfrac{\partial\vec{N}}
{\partial\tau}\right\rangle=0$$ is obtained by
$\dfrac{\partial\vec{N}}{\partial\tau}=0$. Moreover,
$$\aligned S_{\tau\tau}=&\left\langle\dfrac{\partial^{2}\widetilde{F}}{\partial
\tau^{2}},-\vec{N}\right\rangle+\left\langle\dfrac{\partial\widetilde{F}}{\partial
\tau},-\dfrac{\partial\vec{N}}{\partial \tau}\right\rangle\\
=&\left\langle\dfrac{\partial F} {\partial u}\dfrac{\partial^{2}
u}{\partial \tau^{2}}+\dfrac{\partial^{2} F} {\partial
u^{2}}\left(\dfrac{\partial
u}{\partial\tau}\right)^{2}+2\dfrac{\partial^{2}F}{\partial
u\partial t}\dfrac{\partial u}{\partial \tau}+\dfrac{\partial^{2}
F}{\partial t^{2}},-\vec{N}\right\rangle\\
=&\left\langle\dfrac{\partial^{2} F} {\partial
u^{2}}\left(\dfrac{\partial
u}{\partial\tau}\right)^{2}+2\dfrac{\partial^{2}F}{\partial
u\partial t}\dfrac{\partial u}{\partial \tau}+\dfrac{\partial^{2}
F}{\partial t^{2}},-\vec{N}\right\rangle\\
=&\left\langle\dfrac{\partial^{2} F} {\partial
u^{2}}\left(\dfrac{\partial
u}{\partial\tau}\right)^{2}+\dfrac{\partial^{2}F}{\partial u\partial
t}\dfrac{\partial u}{\partial
\tau},-\vec{N}\right\rangle+\left\langle\dfrac{\partial^{2}F}{\partial
u\partial t}\dfrac{\partial u}{\partial \tau}+\dfrac{\partial^{2}
F}{\partial t^{2}},-\vec{N}\right\rangle\\
=&\left\langle\left(\dfrac{\partial F}{\partial
u}\right)_{\tau},-\vec{N}\right\rangle\dfrac{\partial u}{\partial
\tau}+\left \langle\dfrac{\partial^{2}F}{\partial u\partial
t}\dfrac{\partial u}{\partial \tau}+\dfrac{\partial^{2} F}{\partial
t^{2}},-\vec{N}\right\rangle\\=&\left\langle\dfrac{\partial^{2}F}{\partial
u\partial t}\dfrac{\partial u}{\partial
\tau},-\vec{N}\right\rangle-k.\endaligned$$ In terms of the normal
evolving curve, we have $$\left\langle\dfrac{\partial F}{\partial
t},\vec{T}\right\rangle\equiv0\quad {\rm{ for\; all}}\;\;
t\in[0,T).$$ By the formula
$$S_{\tau}=\left\langle\dfrac{\partial F}{\partial
t},-\vec{N}\right\rangle,$$ we get
$$\aligned S_{\theta\tau}=&\left\langle\dfrac{\partial^{2}F}{\partial u\partial t}
\dfrac{\partial u}{\partial \theta},-\vec{N}\right\rangle+
\left\langle\dfrac{\partial F}{\partial
t},-\dfrac{\partial\vec{N}}{\partial\theta}\right\rangle\\
=&\left\langle\dfrac{\partial^{2}F}{\partial u\partial t}
\dfrac{\partial u}{\partial \theta},-\vec{N}\right\rangle+
\left\langle\dfrac{\partial
F}{\partial t},\vec{T}\right\rangle\\
=&\left\langle\dfrac{\partial^{2}F}{\partial u\partial t}
\dfrac{\partial u}{\partial
\theta},-\vec{N}\right\rangle=\dfrac{1}{\partial\theta/\partial
u}\left\langle\dfrac{\partial^{2}F}{\partial u\partial
t},-\vec{N}\right\rangle\\=&\dfrac{1}{(\partial\theta/\partial
s)(\partial s/
\partial u)}\left\langle\dfrac{\partial^{2}F}{\partial u\partial
t},-\vec{N}\right\rangle=\dfrac{1}{k\upsilon}\left\langle\dfrac{\partial^{2}F}{\partial
u\partial t},-\vec{N}\right\rangle.\endaligned$$  Noting that
$$S_{\theta}=\langle\widetilde{F},\vec{T}\rangle,$$ we have
$$\aligned S_{\tau\theta}=&\left\langle\dfrac{\partial \widetilde{F}}{\partial\tau},\vec{T}\right\rangle
+\left\langle\widetilde{F},\dfrac{\partial\vec{T}}{\partial\tau}\right\rangle=\left\langle\dfrac{\partial
\widetilde{F}}{\partial\tau},\vec{T}\right\rangle\\=&\left\langle\dfrac{\partial
F} {\partial u}\dfrac{\partial u}{\partial\tau}+\dfrac{\partial
F}{\partial t},\vec{T}\right\rangle\\
=&\left\langle\dfrac{\partial F} {\partial u}\dfrac{\partial
u}{\partial\tau},\vec{T}\right\rangle=\upsilon\dfrac{\partial
u}{\partial\tau}.\endaligned$$ Hence, the support function $S$
satisfies
\begin{equation}\aligned S_{\tau\tau}=&\left\langle\dfrac{\partial^{2}F}{\partial
u\partial t}\dfrac{\partial u}{\partial
\tau},-\vec{N}\right\rangle-k\\=&k\upsilon\dfrac{\partial
u}{\partial
\tau}S_{\theta\tau}-k\\
=&kS_{\theta\tau}^{2}-k=(S_{\theta\tau}^{2}-1)k,\endaligned\end{equation}
namely,
\begin{equation}S_{\tau\tau}=\dfrac{S_{\theta\tau}^{2}-1}{S_{\theta\theta}+S},\
\ \ \forall\ (\theta,\tau)\in S^{1}\times [0,T).\end{equation} Then,
it follows from $(\ref{1.6})$ that
\begin{equation}\label{2.4}\left\{\aligned &S
S_{\tau\tau}+S_{\tau\tau}S_{\theta\theta}-S_{\theta\tau}^{2}+1=0,\\
&S(\theta, 0)=h(\theta),
\\
&S_{\tau}(\theta,0)=-\widetilde{f}(\theta),\endaligned\right.\end{equation}
where $h$ is the support function of $F_{0}$, and $\widetilde{f} $
is the initial velocity of the initial curve $F_{0}$.

For an unknown function $z=z(\theta,\tau)$ defined for
$(\theta,\tau)\in\mathbb{R}^{2}$, the corresponding {\it
Monge-Amp\`{e}re equation} reads
\begin{equation}\label{2.5}A+Bz_{\tau\tau}+Cz_{\tau\theta}+Dz
_{\theta\theta}
+E(z_{\tau\tau}z_{\theta\theta}-z_{\theta\tau}^{2})=0,\end{equation}
the coefficients $A,\ B,\ C,\ D$ and $E$ depends on $\tau,\ \theta,\
S,\ S_{\tau},\ S_{\theta}$. We say that the equation (\ref{2.5}) is
$\tau$-{\it hyperbolic} for $S$, if
$$\bigtriangleup^{2}(\tau,\theta,z,z_{\tau},z_{\theta})\triangleq C^{2}-4BD+4AE>0$$
and $$ z_{\theta\theta}+B(\tau,\theta,z,z_{\tau},z_{\theta})\neq
0.$$ We state the initial values $z(0,\theta)=z_{0}(\theta),\
z_{\tau}(0,\theta)=z_{1}(\theta)$ for the unknown function on the
$\theta\in [0,2\pi]$. Moreover, we require the following $\tau$-{\it
hyperbolicity} condition:
$$\aligned &\bigtriangleup^{2}(0,\theta,z_{0},z_{1},z^{\prime}_{0})=(C^{2}-4BD+4A)|_{t=0}>0,\\
&z^{\prime\prime}_{0}+B(0,\theta,z_{0},z_{1},z^{\prime}_{0})\neq
0,\endaligned$$ in which $z_{0}^{\prime}=\dfrac{dz_{0}}{d\theta}$
and $z_{0 }^{\prime\prime}=\dfrac{d^{2}z_{0}}{d\theta^{2}}$.

It is easy to see that the equation (\ref{2.4}) is a hyperbolic
Monge-Am\`{e}re equation, in which $$A=1,\ B=S,\ C=D=0,\ E=1.$$ In
fact,
$$\aligned\bigtriangleup^{2}(\tau,\theta,S,S_{\tau},S_{\theta})=&C^{2}-4BD+4A\\
=&0^{2}-4S\times 0+4\times 1=4>0\endaligned$$ and
$$S_{\theta\theta}+B(\tau,\theta,S,S_{\tau},S_{\theta})=S_{\theta\theta}+S=\dfrac{1}{k}\neq
0.$$ Furthermore, if we assume that \ $h(\theta)$ is third and
$\widetilde{f}(\theta)$ is twice continuous by differentiable on the
real axis, then the initial conditions satisfies
$$\bigtriangleup^{2}(0,\theta,h,-\widetilde{f},h_{\theta})=4>0$$ and
$$h_{\theta\theta}+B(0,\theta,h,-\widetilde{f},h_{\theta})=h_{\theta\theta}+h
=\dfrac{1}{k_0}\neq 0.$$
 This implies that the equation (\ref{2.4}) is a
{\it hyperbolic Monge-Amp\`{e}re equation} on $S$. By the standard
theory of hyperbolic equations (e.g., \cite{ha},\ \cite{H},\
\cite{k},\ \cite{L} or \cite{T}), we have

\begin{Theorem}\label{2.1} {\bf (Local existences and uniqueness)}
Suppose that $F_{0}$ is a smooth strictly convex closed curve. Then
there exist a positive $T$ and a family of strictly convex closed
curves $F(\cdot,t)$ (in which $t\in [0,T)$) such that $F(\cdot,t)$
satisfies (\ref{1.6}) (or \ref{2.4}), provided that $f(u)$ is a
smooth function on $S^{1}$.
\end{Theorem}

Theorem 2.1 is nothing but Theorem 1.1, which is one of main results
in this paper.

\section{ An example and some propositions}

In this section, we will give an example to understand further the
normal hyperbolic mean curvature flow. For simplicity we replace
$\tau$ by $t$.

\begin{Example}\label{3.1} Consider $F(\cdot,t)$ to be a family of round circles with the radius R(t)
centered at the origin. The support function and the curvature are
given by $S(\theta,t)=R(t)$ and $k(\theta,t)=1/R(t)$, respectively.
Substituting these into (\ref{1.6}) gives
\begin{equation}\label{3.1}\begin{cases}R_{tt}=-\dfrac{1}{R}~,\\
R(0)=r_0>0,~~~R_{t}(0)=r_1.\end{cases}\end{equation}\end{Example}

For this initial value problem, we have the following lemma which is
given in \cite{h}.

\begin{Lemma}\label{3.1} For arbitrary initial data $r_0>0$,
if the initial velocity $r_1\leqslant 0,$ the solution $R=R(t)$
decreases and shrinks to a point at time $T^{*}$ (where
$T^{*}\leq\sqrt{\frac{\pi}{2}}r_0$, and the equality holds if and
only if $r_{1}=0$); if the initial velocity is positive, the
solution $R$ increases first and then decreases and shrinks to a
point in a finite time.
\end{Lemma}

\begin{Remark}\label{3.1} In fact, this phenomena can also be interpreted by
physical principles. From (\ref{2.4}), we can see that the direction
of acceleration is always the same as the inner normal vector. Thus,
if $R_t(0)\leqslant0$, i.e., the initial velocity is in accordance
with the unit inner normal vector, then evolving circle will shrink
to a point at a finite time; if $r_t(0)>0$, i.e., the initial
velocity is in accordance with the outer unit normal vector, then
the evolving sphere will expand first and then shrink to a point at
a finite time. In the equation $(\ref{2.4})$,\ \
$S_{t}(\theta,0)=-\widetilde{f}\leq0$, i.e., we assume the initial
velocity always accords with the initial unit inner normal, hence
only the first phenomena will happen.\end{Remark}

In what follows, we shall establish some properties enjoyed by the
hyperbolic mean curvature flow.

Consider the following general second-order operator
\begin{equation}\label{3.2}L[w]\triangleq
aw_{\theta\theta}+2bw_{\theta
t}+cw_{tt}+dw_{\theta}+ew_{t},\end{equation} where $a,\ b,\ $ and
$c$ are twice continuously differentiable and $d$ and $e$ are
continuously differentiable functions of $\theta$ and $t$. The
operator $L$ is said to be {\it hyperbolic} at a point $(\theta,t)$,
if
$$b^{2}-ac>0.$$  It is hyperbolic in a domain $D$ if it is
hyperbolic at each point of $D$, and {\it uniformly hyperbolic} in
$D$ if there is a constant $\mu$ such that $$b^{2}-ac\geq\mu>0$$ in
$D$.

We suppose that $w$ and the conormal derivative
$$\dfrac{\partial w}{\partial\nu}\triangleq -b\dfrac{\partial w}
{\partial\theta}-c\dfrac{\partial w}{\partial t}$$ are given at
$t=0$.

We associate with $L$ the {\bf adjoint operator}
$$\aligned L^{*}[\omega]\triangleq&(a\omega)_{\theta\theta}+2(b\omega)_{\theta t}+(c\omega)_{tt}-(d\omega)_{\theta}-(e\omega)_{t}\\
=&a\omega_{\theta\theta}+2b\omega_{\theta
t}+c\omega_{tt}+(2a_{\theta}+2b_{t}-d)\omega_{x}+(2b_{\theta}+2c_{t}-3)\omega_{t}\\
&+(a_{\theta\theta}+2b_{\theta t}+c_{tt}-d_{\theta}-e_{t})\omega.
\endaligned$$

Now we shall show that for any hyperbolic operator $L$ there is a
function $l$ which satisfies the following condition
$$\label{(**)}\left\{\aligned &2\sqrt{b^{2}-ac}\left[l_{t}-\dfrac{1}{c}(\sqrt{b^{2}-ac}-b)l_{\theta}\right]+lK_{+}\geq
0,\\
&2\sqrt{b^{2}-ac}\left[l_{t}+\dfrac{1}{c}(\sqrt{b^{2}-ac}-b)l_{\theta}\right]+lK_{-}\geq
0,\\ & (L^{*}+g)[\omega]\geq 0,\endaligned\right.
$$
 in a sufficiently small strip $0\leq t\leq t_{0}$, where
 $$\aligned K_{+}\triangleq &K_{+}(\theta,t)\triangleq(\sqrt{b^{2}-ac})_{\theta} +\dfrac{b}{c}(\sqrt{b^{2}-ac})_{\theta}+\dfrac{1}{c}(b_{\theta}+c_{t}-e))\sqrt{b^{2}-ac}\\
 &+\left[-\dfrac{1}{2c}(b^{2}-ac)_{\theta}+a_{\theta}+b_{t}-d-\dfrac{b}{c}(b_{\theta}+c_{t}-e)\right],\endaligned$$
 and
$$\aligned K_{-}\triangleq &K_{-}(\theta,t)\triangleq(\sqrt{b^{2}-ac})_{\theta}  +\dfrac{b}{c}(\sqrt{b^{2}-ac})_{\theta}+\dfrac{1}{c}(b_{\theta}+c_{t}-e))\sqrt{b^{2}-ac}\\
 &-\left[-\dfrac{1}{2c}(b^{2}-ac)_{\theta}+a_{\theta}+b_{t}-d-\dfrac{b}{c}(b_{\theta}+c_{t}-e)\right].\endaligned$$ We let
\begin{equation}\label{3.3}l(\theta,t)\triangleq1+\alpha t-\beta t^{2}.\end{equation}

A computation shows that the above condition are
\begin{equation}\label{3.4}\left\{\aligned&2\sqrt{b^{2}-ac}(\alpha-2\beta t)+(1+\alpha
t-\beta t^{2})K_{+}\geq
0,\\
&2\sqrt{b^{2}-ac}(\alpha-2\beta t)+(1+\alpha t-\beta t^{2})K_{-}\geq
0,\\&-2c\beta+(2b_{\theta+2c_{t}-e})(\alpha-2\beta t)\\
&\quad\quad\quad+(a_{\theta\theta}+2b_{\theta
t}+c_{tt}-d_{\theta}-e_{t}+g)(1+\alpha t-\beta t^{2})\geq
0.\endaligned\right.\end{equation} Since all the coefficients and
their derivatives which appear in the above expressions are supposed
bounded and since $-c$ and $\sqrt{b^{2}-ac}$ have positive lower
bounds, the first two expressions above are positive at $t=0$ if
$\alpha$ is chosen sufficiently large. The third expression are
positive at $t=0$ if $\beta$ is chosen sufficiently large. With
these values of $\alpha$ and $\beta$ there is a number $t_{0}>0$
such $l(\theta,t)>0$ and all the inequalities hold for $0\leq t\leq
t_{0}$.

With $l$ given by (\ref{3.3}), the condition on the conormal
derivative becomes
$$\dfrac{\partial\omega}{\partial\nu}+(b_{\theta}+c_{t}-e+c\alpha)\omega\leq 0 \ \ \hbox{at} \ \ t=0.$$
If we select a constant $M$ so large that
\begin{equation}\label{3.5}M\geq -[b_{\theta}+c_{t}-e+ c\alpha]\ \ \
\mbox{on $\;\Gamma_{0}$}.\end{equation} Then we obtain the following
maximum principle for a strip adjacent to the $\theta$-axis(see
\cite{p}).

\begin{Lemma}\label{3.2 }  Suppose that the coefficients of the
operator $L$ given by (\ref{3.2}) are bounded and have bounded first
and second derivatives. Let $D$ be an admissible domain. If $t_{0}$
and $M$ are selected in accordance with (\ref{3.4}) and (\ref{3.5}),
then any function $w$ which satisfies
$$\left\{\aligned&(L+g)[w]\geq 0\ \ \mbox{in D},\\
&\dfrac{\partial w}{\partial\nu}-Mw\leq 0 \ \ \mbox{on
$\Gamma_{0}$},\\
&w\leq 0 \ \ \mbox{on $\Gamma_{0}$}, \endaligned\right.$$ also
satisfies $w\leq 0$ in the part of $D$ which lies in the strip
$0\leq t\leq t_{0}$. The constants $t_{0}$ and $M$ depend only on
lower bounds for $-c$ and $\sqrt{b^{2}-ac}$ and on bounds for the
coefficients of $L$ and their derivatives. \end{Lemma}

\begin{Proposition}\label{3.1}{\bf (Containment principle)} Let $F_{1}$ and $F_{2}: S^{1}\times [0,T)
\rightarrow \mathbb{R}^{2}$ be two convex solutions of $(\ref{1.6})$
(or $(\ref{2.4})$). Suppose that $F_{2}(\cdot,0)$ lies in the domain
enclosed by $F_{1}(\cdot,0)$, and $f_{2}(u)\geq f_{1}(u)$. Then
$F_{2}(\cdot, t)$ is contained in the domain enclosed by
$F_{1}(\cdot,t)$ for all $t\in [0,T)$.
\end{Proposition}
{\bf Proof.} Let $S_{1}(\theta,t)$ and $S_{2}(\theta,t)$ be the
support functions of $F_{1}(\cdot,t)$ and $F_{2}(\cdot,t)$
respectively. Then $S_{1}$ and $S_{2}$ satisfies the same equation
(\ref{2.4}) with $S_{2}(\theta,0)\leq S_{1}(\theta,0)$ and
$S_{2t}(\theta,0)\leq S_{1t}(\theta,0)$ for $\theta\in S^{1}$.

Let $$w(\theta,t)\triangleq S_{2}(\theta,t)- S_{1}(\theta,t),$$ then
$w$ satisfies the following equation
\begin{equation}\label{3.6}\left\{\aligned
&w_{tt}=\left[1-S_{1\theta t} S_{2\theta
t}\right]k_{1}k_{2}w_{\theta\theta}+(k_{1}S_{1\theta
t}+k_{2}S_{2\theta t})w_{\theta t}+\left[1-S_{1}(\theta,t)
S_{2}(\theta,t)\right]k_{1}k_{2}w,\\
&w_{t}(\theta,0)=f_{1}(\theta)-f_{2}(\theta)=w_{1}(\theta),\\
&w(\theta,0)=h_{2}(\theta)-h_{1}(\theta)=w_{0}(\theta).\endaligned\right.\end{equation}
Define the operator $L$ by
\begin{equation}\label{3.7}L[w]\triangleq\left[1-S_{1\theta t}
S_{2\theta t}\right]k_{1}k_{2}w_{\theta\theta}+(k_{1}S_{1\theta
t}+k_{2}S_{2\theta t})w_{\theta t}-w_{tt}.\end{equation} In view of
(3.7), we know that $$a=\left[1-S_{1}(\theta,t)
S_{2}(\theta,t)\right]k_{1}k_{2}, \quad
b=\dfrac{1}{2}(k_{1}S_{1\theta t}+k_{2}S_{2\theta t})\quad
{\rm{and}}\quad c=-1$$ are twice continuously differentiable and
$d=e=0$ are continuously differentiable functions of $\theta$ and
$t$. By a direct computation, we get
$$\aligned b^{2}-ac=&\dfrac{1}{4}(k_{1}S_{1\theta
t}+k_{2}S_{2\theta t})^{2}-\left[1-S_{1}(\theta,t)
S_{2}(\theta,t)\right]k_{1}k_{2}\cdot(-1)\\=&\dfrac{1}{4}(k_{1}S_{1\theta
t}-k_{2}S_{2\theta t})^{2}+k_{1}k_{2}\geq
\min_{\theta\in[0,2\pi]}\{k_{10}(\theta)k_{20}(\theta)\}>0.\endaligned$$
Hence the operator $L$ is defined by (\ref{3.7}) is hyperbolic in
$S^{1}\times [0,T)$ and it is {\it uniformly hyperbolic} in
$S^{1}\times [0,T)$, since there is a constant
$\mu=\min_{\theta\in[0,2\pi]}\{k_{10}(\theta)k_{20}(\theta)\}$ such
that
$b^{2}-ac\geq\mu=\min_{\theta\in[0,2\pi]}\{k_{10}(\theta)k_{20}(\theta)\}>0$
in $S^{1}\times [0,T)$.

By Lemma 3.2, we deduce that
$$S_{2}(\theta,t)\leq S_{1}(\theta,t)$$
for all  $t\in [0,T)$. Thus, the proof is
completed.$\quad\quad\blacksquare$

\begin{Proposition}\label{3.2
} {\bf (Preserving convexity)} Let $k_{0}$ be the mean curvature of
$F_{0}$ and let $\delta=\min_{\theta\in [0,2
\pi]}\{k_{0}(\theta)\}>0$. Then for a $C^{4}$-solution $S$ of
$(\ref{2.4})$, one has $$ k(\theta, t)\geq \delta$$ for $t\in
[0,T_{\max})$, where $[0,T_{\max}) $ is the maximal time interval
for the solution $F(\cdot,t)$ of (\ref{1.6}).\end{Proposition}
\vskip 3mm

\noindent{\bf Proof.} Since the initial curve is strictly convex, by
Theorem 2.1, we know that the solution of $(\ref{2.4})$ remains
strictly convex on some short time interval $[0,T)$ with some $T\leq
T_{\max}$ and its support function satisfies
$$S_{tt}=(S_{\theta t}^{2}-1)k=\dfrac{S_{\theta
t}^{2}-1}{S_{\theta\theta}+S},\quad
 \forall \ (\theta,t)\in S^{1}\times [0,T).$$
By taking derivative in time $t$, we have
$$\aligned
k_{t}=&\left(\dfrac{1}{S_{\theta\theta}+S}\right)_{t}\\=&
-\dfrac{1}{(S_{\theta\theta}+S)^{2}}[S_{\theta\theta t}+S_{t}]\\
=&-k^{2}[S_{\theta\theta t}+S_{t}]\\
=&k^{2}[\widetilde{\sigma}_{\theta\theta}+\widetilde{\sigma}],\endaligned$$
hence $$S_{\theta\theta
t}+S_{t}=-(S_{\theta\theta}+S)^{2}k_{t}=-\dfrac{1}{k^2}k_{t},$$
$$S_{\theta\theta\theta
t}+S_{\theta
t}=\left(-\dfrac{1}{k^{2}}k_{t}\right)_{\theta}=\dfrac{2}{k^{3}}k_{t}k_{\theta}-\dfrac{1}{k^{2}}k_{\theta
t}$$ and
$$\aligned k_{tt}=&\dfrac{2}{(S_{\theta\theta}+S)^{3}}[S_{\theta\theta t}+S_{t}]^{2}-
\dfrac{1}{(S_{\theta\theta}+S)^{2}}[S_{\theta\theta
tt}+S_{tt}]\\=&2k^{3}\left(-\dfrac{1}{k^{2}}k_{t}\right)^{2}-k^{2}\left\{[(S_{\theta
t}^{2}-1)k]_{\theta\theta}+(S_{\theta
t}^{2}-1)k\right\}\\
=&\dfrac{2}{k}k_{t}^{2}-k^{2}\left\{(S_{\theta
t}^{2}-1)_{\theta\theta}k+2(S_{\theta
t}^{2}-1)_{\theta}k_{\theta}+(S_{\theta
t}^{2}-1)k_{\theta\theta}+(S_{\theta
t}^{2}-1)k\right\}\\
=&\dfrac{2}{k}k_{t}^{2}-k^{2}(S_{\theta
t}^{2}-1)(k_{\theta\theta}+k)-k^{2}\left\{2(S_{\theta
t}S_{\theta\theta t})_{\theta}k+4S_{\theta t}S_{\theta\theta
t}k_{\theta}\right\}\\
=&\dfrac{2}{k}k_{t}^{2}-k^{2}(S_{\theta
t}^{2}-1)(k_{\theta\theta}+k)-k^{2}\left\{2(S_{\theta\theta
t}^{2}+S_{\theta t}S_{\theta\theta\theta t})k+4S_{\theta
t}(S_{\theta\theta}+S-S)_{t}k_{\theta}\right\}\\
=&\dfrac{2}{k}k_{t}^{2}-k^{2}(S_{\theta
t}^{2}-1)(k_{\theta\theta}+k)-k^{2}\{2[(S_{\theta\theta
t}+S_{t})^{2}-2S_{\theta\theta t}S_{t}-S_{t}^{2}+S_{\theta
t}(S_{\theta\theta}+S)_{\theta t}-S_{\theta
t}^{2}]k\\
&-4S_{\theta t}\dfrac{1}{k^{2}}k_{t}k_{\theta}-4S_{\theta
t}S_{t}k_{\theta}\}\\
=&\dfrac{2}{k}k_{t}^{2}-k^{2}(S_{\theta
t}^{2}-1)(k_{\theta\theta}+k)-k^{2}\left\{2\left[(S_{\theta\theta
t}+S_{t})^{2}-2(S_{\theta\theta
t}+S_{t})S_{t}+S_{t}^{2}\right.\right.\\
&\left.\left.-S_{\theta t}^{2}+S_{\theta
t}\left(\dfrac{1}{k}\right)_{\theta t}\right]k
 -4S_{\theta t}\dfrac{1}{k^{2}}k_{t}k_{\theta}-4S_{\theta
t}S_{t}k_{\theta}\right\}\\
=&\dfrac{2}{k}k_{t}^{2}-k^{2}(S_{\theta
t}^{2}-1)(k_{\theta\theta}+k)-k^{2}\left\{2\left[\left(-\dfrac{1}{k^{2}}k_{t}\right)^{2}-
2\left(-\dfrac{1}{k^{2}}k_{t}\right)S_{t}+S_{t}^{2}-S_{\theta
t}^{2}\right.\right.\\ &\left.\left.+S_{\theta
t}\left(\dfrac{2}{k^{3}}k_{t}k_{\theta}-\dfrac{1}{k^{2}}k_{\theta
t}\right)\right]k-4S_{\theta
t}\dfrac{1}{k^{2}}k_{t}k_{\theta}-4S_{\theta
t}S_{t}k_{\theta}\right\}\\
=&\dfrac{2}{k}k_{t}^{2}-k^{2}(S_{\theta
t}^{2}-1)(k_{\theta\theta}+k)-k^{2}\left\{\dfrac{2}{k^{3}}k^{2}_{t}+\dfrac{4}{k}S_{t}k_{t}+2S_{t}^{2}-2S_{\theta
t}^{2}k\right.\\&\left.+\dfrac{4}{k^{2}}S_{\theta
t}k_{t}k_{\theta}-\dfrac{2}{k}S_{\theta t}k_{\theta t} -4S_{\theta
t}\dfrac{1}{k^{2}}k_{t}k_{\theta}-4S_{\theta
t}S_{t}k_{\theta}\right\}\\
=&k^{2}(1-S_{\theta t}^{2})k_{\theta\theta}+2kS_{\theta t}k_{\theta
t}+4k^{2}S_{\theta t}S_{t}k_{\theta}-4kS_{t}k_{t}+(S_{\theta
t}^{2}+1-2S_{t}^{2})k^{3}.\endaligned$$ Thus, the curvature $k$
satisfies the following equation
\begin{equation}\label{3.8} k_{tt}=k^{2}(1-S_{\theta t}^{2})k_{\theta\theta}+2kS_{\theta t}k_{\theta
t}+4k^{2}S_{\theta t}S_{t}k_{\theta}-4kS_{t}k_{t}+(S_{\theta
t}^{2}+1-2S_{t}^{2})k^{3}.\end{equation}

Define the operator $L$ as follows
\begin{equation}\label{3.9} L[k]\triangleq k^{2}(1-S_{\theta t}^{2})k_{\theta\theta}+2kS_{\theta t}k_{\theta
t}-k_{tt}+4k^{2}S_{\theta
t}S_{t}k_{\theta}-4kS_{t}k_{t}.\end{equation} In terms of
(\ref{3.6}), $$a=k^{2}(1-S_{\theta t}^{2}), \quad b=kS_{\theta
t}\quad \mbox{and}\quad c=-1$$ are twice continuously differentiable
and $$d=4k^{2}S_{\theta t}S_{t}\quad\mbox{and}\quad e=-4kS_{t}$$ are
continuously differentiable functions of $\theta$ and $t$. By the
direct computation,
$$b^{2}-ac=(kS_{\theta t})^{2}-k^{2}(1-S_{\theta
t}^{2})\cdot (-1)=k^{2}>0,$$ hence the operator $L$ is defined by
(\ref{3.9}) is hyperbolic in the domain $S^{1}\times [0,T)$.

We consider the problem of determining a function $k(\theta,t)$
which satisfies
\begin{equation}\label{3.10}\left\{\aligned &(L+\widetilde{h})[k]
\triangleq  k^{2}(1-S_{\theta t}^{2})k_{\theta\theta}+2kS_{\theta
t}k_{\theta t}+4k^{2}S_{\theta
t}S_{\theta}k_{\theta}\\&\quad\quad\quad\quad\quad-4kS_{t}k_{t}+k^{2}(S_{\theta
t}^{2}+1-2S_{t}^{2})k=0\ \  \mbox{in $S^{1}\times
[0,\widetilde{T})$,}\\
&k(\theta,0)=k_{0}(\theta)\ \ \mbox{on $\Gamma
_{0}$,}\\
&0\leq\dfrac{\partial k}{\partial\nu}\triangleq
-bk_{\theta}-ck_{t}=\gamma(\theta)\ \ \mbox{on $\Gamma
_{0}$.}\endaligned\right.\end{equation} We can find that a function
$\widetilde{k}(\theta, t)=\min_{\theta\in
[0,2\pi]}\left\{k_{0}(\theta)\right\}=\delta$ which satisfies
\begin{equation}\label{3.11}\left\{\aligned
&(L+\widetilde{h})[\widetilde{k}]=0 \ \  \mbox{in $S^{1}\times
[0,\widetilde{T})$,}\\ &\widetilde{k}(\theta,0)\leq k_{0}(\theta)\ \
\mbox{on $\Gamma _{0}$,}\\ &\dfrac{\partial
\widetilde{k}}{\partial\nu}-M\widetilde{k}\leq
\gamma(\theta)-Mk_{0}(\theta)\ \ \mbox{on $\Gamma
_{0}$,}\endaligned\right.\end{equation} where $\Gamma_{0}$ is the
initial domain, and $M$ is the constant given by (\ref{3.5}). If $k$
and if $\widetilde{k}$ satisfies (\ref{3.11}), we may apply Lemma
3.2 to $ \widetilde{k}-k$  and conclude that
$$\widetilde{k}\leq k(\theta,t)\ \ \mbox{in  $S^{1}\times [0,t_{0})$}$$
with $t_{0}\leq T$. This implies that the solution $F(\cdot,t)$ is
convex on $[0, T_{\max})$. Moreover, the curvature of $F(\cdot,t)$
has a uniform positive lower bound $\min_{\theta\in S^{1}}\{
k_{0}(\theta)\}$ on $S^{1}\times [0,T_{\max})$. Thus, the proof is
completed.$\quad\quad\blacksquare$

The following lemmas will be useful later.
\begin{Lemma}\label{3.3} The arclength $\mathscr{L}(t)$ of the closed
curve $F(\cdot,t)$ satisfies
$$\dfrac{d \mathscr{L}(t)}{
dt}=-\int_0^{2\pi}\widetilde{\sigma}(\theta,t)d\theta$$ and
$$\dfrac{d^{2}\mathscr{L}(t)}{dt^{2}}=\int_0^{2\pi}\left[
\left(\dfrac{\partial\widetilde{\sigma}}{\partial\theta}\right)^{2}k-k\right]d\theta.$$\end{Lemma}
{\bf Proof.} By the definitino of arclength,
$$\mathscr{L}(t)=\int_0^{2\pi}v(\theta,t)d\theta.$$
By a direct calculation,
$$\aligned\dfrac{d \mathscr{L}(t)}{
dt}=&\int_0^{2\pi}\dfrac{\partial v}{\partial
t}d\theta\\=&\int_0^{2\pi}-\widetilde{\sigma}(\theta,t)k(\theta,t)v(\theta,t)d\theta\\
=&-\int_0^{2\pi}\widetilde{\sigma}(\theta,t)d\theta,\endaligned$$
and then
$$\aligned\dfrac{d^{2}\mathscr{L}(t)}{dt^{2}}=&-\int_0^{2\pi}\dfrac{\partial }{\partial t}
\left(\widetilde{\sigma}(\theta,t)\right)d\theta\\
=&\int_0^{2\pi} \left(S_{\theta
t}^{2}-1\right)kd\theta\\
=&\int_0^{2\pi}\left[
\left(\dfrac{\partial\widetilde{\sigma}}{\partial\theta}\right)^{2}k-k\right]d\theta.\endaligned$$
Thus, the proof is completed.$\quad\quad\blacksquare$

\begin{Lemma}\label{3.4} The area $\mathscr{A}(t)$ enclosed by the closed curve $F(\cdot,t)$ satisfies
 $$\aligned&\dfrac{d
\mathscr{A}(t)}{dt}=\int_0^{2\pi}\dfrac{S_{t}}{k}d\theta,\\
&\dfrac{d^{2}
\mathscr{A}(t)}{dt^{2}}=-2\pi+\int_0^{2\pi}S_{t}^{2}d\theta,\\
&\dfrac{d^{3} \mathscr{A}(t)}{dt^{3}}=\int_0^{2\pi}(S_{\theta
t}^{2}-1)kS_{t}d\theta.\endaligned$$\end{Lemma}

\noindent{\bf Proof.} The area $\mathscr{A}(t)$ enclosed by the
convex curve is defined by
$$\aligned \mathscr{A}(t)=&-\dfrac{1}{2}\int_0^{2\pi}\left<\widetilde{F}(\theta,t),v(\theta,t)N(\theta,t)\right>d\theta\\
=&\dfrac{1}{2}\int_0^{2\pi}\dfrac{S}{k}d\theta.\endaligned$$ Then,
$$\aligned\dfrac{d
\mathscr{A}(t)}{dt}=&\dfrac{1}{2}\int_0^{2\pi}\left[\dfrac{S_{t}}{k}-\dfrac{S}{k^{2}}k_{t}\right]d\theta,\\
=&\dfrac{1}{2}\int_0^{2\pi}\left[\dfrac{S_{t}}{k}+S(S_{\theta\theta
t}+S_{t})\right]d\theta\\
=&\dfrac{1}{2}\int_0^{2\pi}\left[\dfrac{S_{t}}{k}+(S_{\theta\theta
}+S)S_{t}\right]d\theta\\
=&\dfrac{1}{2}\int_0^{2\pi}\left[\dfrac{S_{t}}{k}+\dfrac{S_{t}}{k}\right]d\theta\\
=&\int_0^{2\pi}\dfrac{S_{t}}{k}d\theta,\endaligned$$ and then,
$$\aligned\dfrac{d^{2}
\mathscr{A}(t)}{dt^{2}}=&\int_0^{2\pi}\dfrac{\partial}{\partial
t}\left(\dfrac{S_{t}}{k}d\theta\right)\\
=&\int_0^{2\pi}\left[\dfrac{S_{tt}}{k}-\dfrac{S_{t}}{k^{2}}k_{t}\right]d\theta\\
=&\int_0^{2\pi}\left[S_{\theta t}^{2}-1+S_{t}(S_{\theta\theta
t}+S_{t})\right]d\theta\\
=&\int_0^{2\pi}(S_{\theta t}^{2}-1+S_{t}^{2}-S_{\theta
t}^{2})d\theta\\
=&-2\pi+\int_0^{2\pi}S_{t}^{2}d\theta.\endaligned$$ Finally,
$$\aligned\dfrac{d^{3}
\mathscr{A}(t)}{dt^{3}}=&2\int_0^{2\pi}S_{t}S_{tt}d\theta\\
=&2\int_0^{2\pi}S_{t}(S_{\theta t}^{2}-1)kd\theta.\endaligned$$
Thus, the proof is completed.$\quad\quad\blacksquare$

\begin{Lemma}\label{3.5} Under the Proposition 3.2, the following
inequality holds
$$\left(\dfrac{\partial\widetilde{\sigma}}{\partial\theta}\right)^{2}-1<0
\quad {\rm{for\ all}}\quad t\in [0,T_{\max}).$$\end{Lemma}

\noindent{\bf Proof.} Since
$$\dfrac{\partial\sigma}{\partial t}=k>0\quad {\rm{for\ \ all}}\quad
t\in [0,T_{\max}),$$ then
$$\sigma(u,t)>\sigma(u,0)\quad {\rm{for\ all}}\ \ t\in (0,T_{\max}),$$
i.e.,
$$\widetilde{\sigma}(\theta,t)=\sigma(u,t)>\sigma(u,0)=\widetilde{\sigma}(\theta,0)\quad
{\rm{for\ \ all}}\quad t\in (0,T_{\max}),$$ hence, $$\dfrac{\partial
\widetilde{\sigma}}{\partial t}>0\quad {\rm{for \ all}}\quad t\in
[0,T_{\max}).$$On the other hand, by the chain rule,
$$\aligned\dfrac{\partial
\sigma}{\partial t}=&\dfrac{\partial\widetilde{\sigma}}
{\partial\theta}\dfrac{\partial\theta}{\partial
t}+\dfrac{\partial\widetilde{\sigma}}{\partial
t}\\=&\dfrac{\partial\widetilde{\sigma}}
{\partial\theta}\dfrac{\partial\sigma}{\partial
s}+\dfrac{\partial\widetilde{\sigma}}{\partial
t}\\=&\dfrac{\partial\widetilde{\sigma}}
{\partial\theta}\dfrac{\partial\widetilde{\sigma}}
{\partial\theta}\dfrac{\partial\theta}{\partial
s}+\dfrac{\partial\widetilde{\sigma}}{\partial t},\endaligned$$
hence,  $$\dfrac{\partial\widetilde{\sigma}}{\partial t}
=\left[1-\left(\dfrac{\partial\widetilde{\sigma}}
{\partial\theta}\right)^{2}\right]k> 0,$$ therefore,
$$\left(\dfrac{\partial\widetilde{\sigma}} {\partial\theta}\right)^{2}-1< 0 \quad {\rm{
for\ \ all}}\quad t\in [0,T_{\max}).$$ Thus, the proof is
completed.$\quad\quad\blacksquare$

\section{Shrinking to a point --- Proof of Theorem 1.2}

From Example 3.1, we know that, when $t\to T^{*}$, the solution
$F(\cdot,t)$ converges to a point. In this section, we will show
that such phenomenon actually holds for all the evolutions of
strictly convex closed curves with suitably initial velocities. In
other words, we will prove the following theorem which implies
Theorem 1.2.

\begin{Theorem}\label{4.1} Suppose that $F_{0}$ is a smooth strictly
convex closed curve and the initial velocity $f$ satisfies all
assumptions mentioned in Section 3. Then the solution of (\ref{1.6})
with $F_{0}$ and $f$ as initial curve and initial velocity of the
initial curve, respectively, exists only at a finite time interval
$[0,T_{\max})$. Moreover, when $t\to T_{\max}$, the solution
$F(\cdot,t)$ converges to a point.\end{Theorem}

 \noindent{\bf Proof.} Let $[0,T_{\max})$ be the
maximal time interval for the solution $F(\cdot,t)$ of (\ref{1.6})
with $F_{0}$ and $f$ as initial curve and initial velocity of the
initial curve, respectively. We divide the proof into four steps.

\vskip 3mm

\noindent{\bf Step 1. Preserving convexity}

By Proposition 3.2, we know that the solution $F(\cdot,t)$ remains
strictly convex on $[0,T_{\max})$ and the curvature of $F(\cdot,t)$
has a uniform positive lower bound $\min_{\theta\in
S^{1}}\{k_{0}(\theta)\}$ on $S^{1}\times [0,T_{\max})$.

\vskip 3mm

\noindent{\bf Step 2. Finite time existence}

Enclose the initial curve $F_{0}$ by a large circle $\gamma_{0}$
with the normal initial velocity equals to the normal initial
velocity of the initial curve $F_{0}$. Then evolve $\gamma_{0}$ by
the flow $(\ref{1.6})$ to get a solution $\gamma(\cdot,t)$. By
Example 3.1, we know that the solution $\gamma(\cdot,t)$ exists only
at a finite time interval $[0,T^{*})$, and $\gamma(\cdot,t)$
converges to a point when $t\to T^{*}(<+\infty)$. Applying
Proposition 3.1 (containment principle), we deduce that $F(\cdot,t)$
is always enclosed by $\gamma(\cdot,t)$ for all $t\in[0,T^{*})$.
Thus we conclude that the solution $F(\cdot,t)$ must become singular
at some time $T_{\max}\leq T^{*}$.

\vskip 3mm

\noindent{\bf Step 3. Hausdorff convergence}

Note that $F(\cdot, t_{2})$ is enclosed by $F(\cdot,t)$ whenever
$t_{2}>t_{1}$ by the evolution equation $(\ref{1.6})$ (or
(\ref{2.4})). In other words, $F(\cdot,t)$ is shrinking. Let us
recall the following classical result in convex geometry (see
\cite{s}).

\noindent{\bf Blaschke Selection Theorem}{\it \ Let $K_{j}$ be a
sequence of convex sets which are contained in a bounded set. Then
there exists a subsequence ${K_{j_{k}}}$ and a convex set $K$ such
that ${K_{j_{k}}}$ converges to $K$ in the Hausdorff metric.}

Thus, by using this result, we can directly deduce that $F(\cdot,t)$
converges to a (maybe degenerate and nonsmooth) weakly convex curve
$F(\cdot,T_{\max})$ in the Haufdorff metric.

\vskip 3mm

\noindent{\bf Step 4. Shrinking to a point}

Noting that
$$\left(\dfrac{\partial\tilde{\sigma}}{\partial\theta}\right)^{2}-1<0\quad {\rm{ for\; all}}\quad t\in [0,T_{\max}),$$
we obtain from Lemma 3.2 that
$$\dfrac{d^{2}\mathscr{L}(t)}{dt^{2}}<0,\ \quad \dfrac{d\mathscr{L}(t)}{dt}<0
\quad {\rm{ for\; all}}\quad t\in [0,T_{\max}).$$ Hence, there
exists a finite time $T_{0}$ such that $\mathscr{L}(T_{0})=0$,
provided that $T_{0}\le T_{\max}$. There will be two cases:

\vskip 2mm

\noindent {\bf Case I:} $T_{0}\leq T_{\max}$.  On the one hand,
there exists a unique classical solution of the Cauchy problem
(\ref{1.6}) on the interval $[0,T_0)$; on the other hand, when $t$
goes to $T_0$, $\mathscr{L}(t)$ tends to zero, i.e.,
$$\mathscr{L}(t)\longrightarrow 0\quad {as}\quad t\nearrow T_0.$$
This implies that the curvature $k$ goes to infinity when $t$ tends
to $T_0$, and then the solution will blow up at the time $T_0$.
Therefore, by the definition of $T_{\max}$, we have
$$T_{0}=T_{\max}.$$ That is, when $t\nearrow T_{\max}$, the solution
$F(\cdot,t)$ converges to a point.

\vskip 2mm

\noindent {\bf Case I\!I:} $T_{0}>T_{\max}$. In the present
situation, $$\mathscr{L}(T_{\max})> 0.$$ Then $F(\cdot,T_{\max})$
must be a line segment. It is clear that $\min_{\theta\in
S^{1}}\{k(\theta,t)\}$ tends to zero, when $t$ goes to $T_{max}$.
But in Step 1, we have shown that the curvature of $F(\cdot,t)$ has
a uniform positive lower bound. Hence, Case I\!I is not possible.
Thus, the proof of Theorem 4.1 is completed.$\quad\quad\blacksquare$

\section{Normal hyperbolic mean curvature flow with dissipation}
In this section, we consider the normal hyperbolic mean curvature
flow with dissipative term
\begin{equation}\label{5.1}\left\{\aligned&\dfrac{\partial^{2} F}{\partial t^{2}}(u,t)
=k(u,t)\vec{N}(u,t)-\triangledown \rho(u,t)+d\dfrac{\partial F}{\partial t},\quad \forall \ (u,t)\in S^{1}\times[0,T),\\
&F(u,0)= F_{0}(u),\\
&\dfrac{\partial F}{\partial t}(u,0)=
f(u)\vec{N_0},\endaligned\right.\end{equation} where $k$ is the mean
curvature, $\vec{N}$ is the inner unit normal at $F(u,t)$, $F_{0}$
stands for the initial strictly convex smooth closed curve, $f(u)$
and $\vec{N_{0}}$ are the initial velocity and inner normal vector
of $F_{0}$, respectively, $d$ is a negative constant and
$\triangledown\rho$ is denoted by
$$\triangledown \rho\triangleq\left\langle\dfrac{\partial^{2} F}{\partial s\partial t},
\dfrac{\partial F}{\partial t}\right\rangle\vec{T}.$$

For $(\ref{5.1})$, if we assume $$\left\langle\dfrac{\partial
F}{\partial t},\vec{N}\right\rangle=v,$$ then we obtain from
$(\ref{5.1})$ that
\begin{equation}\label{5.2}v_{t}=k+dv,\end{equation}which is the same as the equation $(\ref{1.2})$.

Similar to Section 2, we can also derive a hyperbolic
Monge-Amp\`{e}re equation.

In fact, let us use the normal angle to parameterize each convex
curve $F(\cdot,t)$, that is, set
$$\widetilde{F}(\theta,\tau)=F(u(\theta,\tau), t(\theta,\tau)),$$
where $t(\theta,\tau)=\tau$. Here $\vec{N},\ \vec{T}$ and $\theta$
are independent of the parameter $\tau$. Let $$\dfrac{\partial
F}{\partial t}=a(u,t)\vec{N}+b(u,t)\vec{T}\quad {\rm{for\ \
all}}\quad t\in [0,T),$$ then we have
$$\aligned \dfrac{\partial\theta}{\partial t}=&\dfrac{\partial a(u,t)}{\partial
s}+k(u,t)b(u,t)\\=&\dfrac{\partial
\widetilde{a}(\theta,\tau)}{\partial \theta}\dfrac{\partial
\theta}{\partial s}+k(\theta,\tau)\widetilde{b}(\theta,\tau)\\
=&k(\theta,\tau)\left(\dfrac{\partial
\widetilde{a}(\theta,\tau)}{\partial
\theta}+\widetilde{b}(\theta,\tau)\right).\endaligned$$ On the other
hand, the support function satisfies $$S_{\tau}=-\widetilde{a},$$
hence,
$$\dfrac{\partial\theta}{\partial
t}=k\left(\dfrac{\partial}{\partial\theta}\left(-S_{\tau}\right)+\widetilde{b}\right
)=k(-S_{\theta\tau}+\widetilde{b}).$$ Then the support function
$S(\theta,\tau)$ satisfies the following equation
$$S_{\tau\tau}=[(S_{\theta\tau}-\widetilde{b})^{2}-1]k+dS_{\tau},$$
namely,
\begin{equation}\label{5.3} S_{\tau\tau}=\dfrac{(S_{\theta\tau}-\widetilde{b})^{2}-1}{S_{\theta\theta}+S}+dS_{\tau}. \end{equation}
The equation $(\ref{5.3})$ is equivalent to the following equation
\begin{equation}\label{5.4} S S_{\tau\tau}+2\widetilde{b}S_{\theta\tau}-dS_{\tau}S_{\theta\theta}+
S_{\tau\tau}S_{\theta\theta}-S_{\theta\tau}^{2}+1-\widetilde{b}^{2}-dSS_{\tau}=0.\end{equation}
Denote $$A=1-\widetilde{b}^{2}-dSS_{\tau},\quad B=S,\quad
C=2\widetilde{b},\quad D=-dS_{\tau}\quad{\rm{and}}\quad E=1,$$ then
$$\aligned \bigtriangleup^{2}(\tau,\theta,S,S_{\tau},S_{\theta})=&C^{2}-4BD+4AE
\\=&(2\widetilde{b})^{2}
-4S(-dS_{\tau})+4(1-\widetilde{b}^{2}-dSS_{\tau})\\=&4>0,\\
S_{\theta\theta}+B(\tau,\theta,S,S_{\tau},S_{\theta})=&S_{\theta\theta}+S=\dfrac{1}{k}\neq
0.\endaligned$$ Furthermore, we state the initial values
$S(0,\theta)=h_{0}(\theta),\
S_{\tau}(0,\theta)=-\widetilde{f}(\theta)$ for the unknown function
on the $\theta\in [0,2\pi]$,\ $h(\theta)$ being third and
$\widetilde{f}(\theta)$ twice continuous by differentiable on the
real axis. Moreover, we require the  $\tau$-hyperbolicity condition:
$$\aligned &\bigtriangleup^{2}(0,\theta,\epsilon h,-\epsilon\widetilde{f},h_{\theta})=(C^{2}-4BD+4A)|_{t=0}=4>0,\\
&S_{0\theta\theta}+B(0,\theta,\epsilon h
,-\epsilon\widetilde{f},h_{\theta})=h_{\theta\theta}+h=\dfrac{1}{k_{0}}\neq
0.\endaligned$$ This implies the equation $(\ref{5.4})$ is a
hyperbolic Monge-Amp\`{e}re equation on $S$. Then the support
function $S$ satisfies the following initial value problem
\begin{equation}\label{5.5}\left\{\aligned &S
S_{\tau\tau}+2\widetilde{b}S_{\theta\tau}-dS_{\theta\tau}+S_{\tau\tau}S_{\theta\theta}
-S_{\theta\tau}^{2}+1-\widetilde{b}^{2}-dS_{\tau}=0,\\
&S(\theta, 0)= h(\theta),\\
&S_{\tau}(\theta,0)=-\widetilde{f}(\theta),\endaligned\right.\end{equation}
where $h$ is the support function of $F_{0}$, and $\widetilde{f} $
is the initial velocity of the initial curve $F_{0}$.

Similarly, we can get the curvature the curvature $k$ satisfies the
following equation,
\begin{equation}\label{5.5}\aligned k_{\tau\tau}=&k^{2}[1-(S_{\theta\tau}-\tilde{b})^{2}])k_{\theta\theta}
+2k(S_{\theta\tau}-\tilde{b})k_{\theta\tau} +4k^{2}(S_{\theta
\tau}-\tilde{b})(S_{\tau}+\tilde{b}_{\theta}) k_{\theta}\\&
+(d-4kS_{\tau}-4k\tilde{b}_{\theta})k_{\tau}+[S_{\theta
\tau}^{2}+1-2S_{\tau}^{2}-4s_{\tau}\tilde{b}_{\theta}-\tilde{b}^{2}
+2(S_{\theta\tau}-\tilde{b})\tilde{b}_{\theta\theta}]k^{3}.\endaligned\end{equation}
It is easy to verify that the equation (\ref{5.5}) is also a
hyperbolic equation about $\tau$.

Motivated by the theory of dissipative hyperbolic equations (see
\cite{k0},\ \cite{N}), we will study the initial value problem (5.5)
and the initial value problem for  (5.6) in the forthcoming paper.

\section{Relation between the hyperbolic mean curvature flow and the string evolution equation
 in the Minkowski space $\mathbb{R}^{1,1}$}

In this section, we study the relation between the hyperbolic mean
curve flow and the evolution equation for the string in the
Minkowski space $\mathbb{R}^{1,1}$.

Let $z=(z_{0},z_{1})$ be a position vector of a point in the
two-dimensional Minkowski space $\mathbb{R}^{1,1}$. The scalar
product of two vectors $z$ and $w=(w_{0},w_{1})$ is
$$\left<z,w\right>=-z_{0}w_{0}+z_{1}w_{1}.$$
The Lorentz metric of $\mathbb{R}^{1,1}$ reads
$$ds^{2}=-dt^{2}+du^{2}.$$

A massless closed curve moving in two-dimensional Minkowski space
can be defined by making its action proportional to the
two-dimensional area swept out in the Minkowski space. We are
interested in the following motion of one-dimensional Riemannian
manifold in $\mathbb{R}^{1,2}$ with the following parameter
\begin{equation}\label{6.1}(t,u)\to \widetilde{X}=(t,X(t,u)),
\end{equation} where $u\in \mathscr{M}$ and $\widetilde{X}(\cdot,t)$ be a positive
vector of a point in the Minkowski space $\mathbb{R}^{1,2}$. The
induced Lorentz metric reads
\begin{equation}\label{6.2}\left\{\aligned\widetilde{g}_{00}=&-1+\left(\dfrac{\partial X}{\partial t},\dfrac{\partial X}{\partial t}\right),\\
\widetilde{g}_{01}=&\widetilde{g}_{10}=\left(\dfrac{\partial
X}{\partial t},\dfrac{\partial X}{\partial
u}\right),\\
\widetilde{g}_{11}=&g_{11}=\left(\dfrac{\partial X}{\partial
u},\dfrac{\partial X}{\partial u}\right),\endaligned\right.
\end{equation}
i.e., the Lorentz metric becomes
$$ds^{2}=(dt,du)A(dt,du)^{T},$$
where \begin{equation}A=\left(\begin{array}{cc}|X_t|^{2}-1
 &\langle X_t,X_u\rangle\\
\langle X_t,X_u\rangle\ \ \ \ \
&|X_{u}|^{2}\end{array}\right),\end{equation} in which
$$|X_{t}|^{2}=\langle X_{t},X_{t}\rangle,\ \ \ |X_{u}|^{2}=\langle X_{u},X_{u}\rangle.$$
Kong, Zhang and Zhou in \cite{k1} investigated the dynamics of
relativistic (in particular, closed) strings moving in the Minkowski
space $\mathbb{R}^{1,n}(n\geq 2)$. By the variational method, they
get the following equation
\begin{equation}\label{6.3}|X_{u}|^{2}X_{tt}-2\langle X_{t},X_{u}\rangle X_{tu}+(|X_{t}|^{2}-1)X_{uu}=0.\end{equation}
 Except the variational method, by vanishing mean curvature of the
sub-manifold $\mathscr{M}$, we can obtain the following equation for
the motion of $\mathscr{M}$ in the Minkowski space
$\mathbb{R}^{1,2}$
\begin{equation}\label{6.4}\widetilde{g}^{\alpha\beta}\triangledown_{\alpha}\triangledown_{\beta}\widetilde{X}=
\widetilde{g}^{\alpha\beta}\left(\dfrac{\partial^{2}\widetilde{X}}{\partial
x^{\alpha}\partial
x^{\beta}}-\widetilde{\Gamma}_{\alpha\beta}^{\gamma}\dfrac{\partial\widetilde{X}}{\partial
x^{\gamma}}\right)=0,\end{equation} where $\alpha,\beta=0,1$. It is
convenient to fix the parametrization partially (see  Albrecht and
Turok \cite{a}, Turok and Bhattacharjee \cite{Tu}) by requiring
\begin{equation}\label{6.5}\widetilde{g}_{01}=\widetilde{g}_{10}=\left(\dfrac{\partial
X}{\partial t},\dfrac{\partial X}{\partial
u}\right)=0,\end{equation} that is, we require the additional gauge
condition that the string velocity be orthogonal to the string
tangent direction. We assume that the surface is $C^{2}$ and {\it
time-like}, i.e.,
$$(|X_{t}|^{2}-1)|X_{u}|^{2}-\langle X_{t},X_{u}\rangle^{2}<0,$$
 equivalently,
$$(1-|X_{t}|^{2})>0.$$
Obviously, the equation $(\ref{6.4})$ is equivalent to
\begin{equation}\label{6.6}\left(\dfrac{\partial^{2}X}{\partial t^{2}},\dfrac{\partial X}{\partial t}\right)
-g^{11}\left(|X_{t}|^{2}-1\right)\left(\dfrac{\partial^{2}X}{\partial
t\partial u},\dfrac{\partial X}{\partial u}\right)=0,\end{equation}
\begin{equation}\label{6.7}\aligned
\dfrac{\partial^{2} X}{\partial
t^{2}}+&g^{11}\left(\dfrac{\partial^{2} X} {\partial
u^{2}}-\Gamma_{11}^{1}\dfrac{\partial X}{\partial
u}\right)\left(|X_{t}|^{2}-1\right)-\dfrac{1}{|X_{t}|^{2}-1}\left(\dfrac{\partial^{2}X}{\partial
t^{2}},\dfrac{\partial X}{\partial t}\right)\dfrac{\partial
X}{\partial t}\\ \ \ \ \
&+g^{11}\left(\dfrac{\partial^{2}X}{\partial t\partial
u},\dfrac{\partial X}{\partial t}\right)\dfrac{\partial X}{\partial
u}+g^{11}\left(\dfrac{\partial^{2}X}{\partial t\partial
u},\dfrac{\partial X}{\partial u}\right)\dfrac{\partial X}{\partial
t}=0.\endaligned
\end{equation}
It is easy to verify that the system $(\ref{6.6})-(\ref{6.7})$ is
equivalent to $(\ref{6.6})$ and the following equation
$$\dfrac{\partial^{2} X}{\partial
t^{2}}+g^{11}\left(\dfrac{\partial^{2} X} {\partial
u^{2}}-\Gamma_{11}^{1}\dfrac{\partial X}{\partial
u}\right)\left(|X_{t}|^{2}-1\right)+g^{11}\left(\dfrac{\partial^{2}X}{\partial
t\partial u},\dfrac{\partial X}{\partial t}\right)\dfrac{\partial
X}{\partial u}=0,$$ namely,
 \begin{equation}\label{6.9}
 \dfrac{\partial^{2}
X}{\partial
t^{2}}=(1-|X_{t}|^{2})k\vec{N}-\dfrac{1}{|X_{u}|^{2}}\left(\dfrac{\partial^{2}X}{\partial
t\partial u},\dfrac{\partial X}{\partial t}\right)\dfrac{\partial
X}{\partial u}.\end{equation}

\begin{Remark}\label{6.1} The equation $(\ref{6.9})$ is similar to the equation in
$(\ref{1.6})$, both of them evolve normally.
 The difference between the equation in $(\ref{1.6})$ and the equation
$(\ref{6.9})$ is only the normal acceleration of the evolving curve.
Because $$1-|X_{t}|^{2}>0,$$ that is, the velocity of the string is
always less than the velocity of light which is meaning in the
classical physics, the motion of the string in the Minkowski space
$\mathbb{R}^{1,1}$ can be regarded as one of applications of general
normal hyperbolic mean curvature flow.\end{Remark}

\vskip 5mm\noindent{\Large {\bf Acknowledgements.}}
 Wang would
like to thank the Center of Mathematical Sciences at Zhejiang
University for the great support and hospitality. The work of Kong
and Wang was supported in part by the NNSF of China (Grant No.
10671124) and the NCET of China (Grant No. NCET-05-0390); the work
of Liu was supported in part by the NSF and NSF of China.

\end{document}